\documentclass[letterpaper, 10 pt, conference]{ieeeconf} 	
\IEEEoverridecommandlockouts    % This command is only needed if 
\newtheorem{theorem}{Theorem}[section]
\newtheorem{lemma}[theorem]{Lemma}
\newtheorem{definition}[theorem]{Definition}

\newtheorem{proposition}[theorem]{Proposition}

\newtheorem{remark}[theorem]{Remark}
\newtheorem{assumption}{Assumption}
\newtheorem{example}[theorem]{Example}
\newenvironment{pfof}[1]{\vspace{1ex}\noindent{\itshape 
	Proof of #1:}\hspace{0.5em}} {\hfill\QEDBL\vspace{1ex}}

\newcommand{\mc}{\mathcal}

\newcommand{\real}{\mathbb{R}} 

\newcommand{\naturalpos}{\mathbb{N}_{>0}}
\newcommand{\realpos}{\mathbb{R}_{> 0}}
\newcommand{\realnneg}{\mathbb{R}_{\geq 0}}

\newcommand{\tsp}{\mathsf{T}} 
 
\newcommand{\inv}{{\negat 1}} 
\newcommand{\negat}{\scalebox{0.75}[.9]{\( - \)}}
\newcommand*{\QEDB}{\hfill\ensuremath{\square}}%  empty square
\newcommand*{\QEDBL}{\hfill\ensuremath{\blacksquare}}% Black square
\newcommand{\diag}{\operatorname{diag}}

\renewcommand{\ast}{\text{\ding{91}}}

\newcommand{\map}[3]{#1: #2 \rightarrow #3}

\newcommand{\sbs}[2]{{#1}_{\textup{#2}}}
\newcommand{\sps}[2]{{#1}^{\textup{#2}}}
\newcommand{\until}[1]{\{1,\dots,#1\}}
\newcommand{\norm}[1]{\Vert #1 \Vert}

\DeclareMathAlphabet{\mymathbb}{U}{BOONDOX-ds}{m}{n}

\usepackage[figurename=Fig.,font=footnotesize]{caption}
\newcommand{\inprod}[2]{\langle  #1 , #2 \rangle}

\usepackage{hyperref}
\usepackage{graphicx,color}
\graphicspath{{./IMG/}{images/},{./IMG-BIOS//}}
\usepackage{amsmath}
\usepackage{amssymb}
\usepackage{mathrsfs}
\usepackage{subfigure}
\usepackage{url}
\usepackage{booktabs}
\usepackage{array}
\usepackage[table]{xcolor}
\usepackage{mathtools} 		% needed for cases environment
\usepackage{epstopdf}
\usepackage{cite}
\usepackage{algorithm}
\usepackage{algorithmic}
\usepackage{upgreek}
\usepackage{dsfont}
\usepackage{pifont} % Needed for \ding

\begin{document}

%%%%%%%%%%%%%%%%%%%%%%%%%%%%%%%%%%%%%%%%%%
\title{\LARGE \textbf{Navigation Systems May Deteriorate Stability in Traffic Networks}}

\author{Gianluca Bianchin and Fabio Pasqualetti}
%%%%%%%%%%%%%%%%%%%%%%%%%%%%%%%%%%%%%%%%%%

\maketitle

\begin{abstract}
Advanced traffic navigation systems, which provide routing recommendations to 
drivers based on real-time congestion information, are nowadays widely adopted 
by roadway transportation users. 
Yet, the emerging effects on the traffic dynamics originating from the widespread 
adoption of these tools have remained largely unexplored until now. 
In this paper, we propose a dynamic model where drivers imitate the path 
preferences of previous drivers, and we study the properties of its 
equilibrium points. Our model is a dynamic generalization of the classical 
\textit{traffic assignment framework}, and extends it by accounting for dynamics 
both in the path decision process and in the network's traffic flows.
We show that when travelers learn shortest paths by imitating other travelers, 
the overall traffic system benefits from this mechanism and transfers the maximum 
admissible amount of traffic demand. 
On the other hand, we demonstrate that when the travel delay functions are not 
sufficiently steep or the rates at which drivers imitate previous travelers are 
not adequately chosen, the trajectories of the traffic system may fail to 
converge to an equilibrium point, thus failing asymptotic stability. 
Illustrative numerical simulations combined with empirical data from highway 
sensors illustrate our findings.
\end{abstract}

\section{INTRODUCTION}

Roadway traffic networks are fundamental components of modern societies, making 
economic activity possible by enabling the transfer of passengers, goods, and 
services in a timely and reliable fashion. 
Despite their critical role, these transportation systems 
are impaired by the long-standing problem of traffic congestion, which 
wastes billions of gallons of fuel each year~\cite{MJS-MC-SSP:18,EC:20}.
Advanced navigation systems are nowadays widely adopted by travelers, largely 
thanks to the widespread use of smartphone-based navigation apps (such as 
Google Maps, Inrix, Waze, Apple Maps, etc.) \cite{JH-DW-RH-XB-QJ-AB:10}. 
Advanced navigation systems provide shortest-path routing recommendations based 
on real-time global travel time information. 
On the one hand, these technologies have enabled travelers to save 
time and fuel but, on the other hand, they have transformed the transportation 
infrastructure originating unanticipated effects and disrupting existing traffic 
flow patterns~\cite{AK-NL-FF-HS-VC-AB-KJ:19}.
While the implications of the widespread adoption of advanced navigation 
systems have been analyzed game-theoretically~\cite{MP:15}, a characterization of 
the impact of these technologies on dynamic models of traffic for general, 
dynamic, traffic networks has remained elusive until now.

In this work, we study the stability properties of a traffic system composed of 
the interconnection between a dynamic model of traffic flows (inspired from 
the Cell Transmission Model~\cite{CFD:95}) and a dynamic model of route 
selection (derived from the Replicator Dynamics~\cite{JWW:97}).
Our choice of using the replicator equation is motivated by recent studies that 
showed that this model emerges as an aggregate description of learning 
processes in large populations and as the limiting case of the best response 
dynamics~\cite{TB-RS:97}. We show that, at equilibrium, our model shares the 
same properties as the well-studied routing game~\cite{MP:15}, and thus it is 
consistent with existing studies that focus on systems operating at 
equilibrium. 
It is worth noting that, with respect to the classical routing-game framework, 
our model accounts for dynamics both in the \textit{route selection process} as 
well as in the \textit{traffic flow model.}
Our dynamical model suggests that systems where travelers continuously prefer 
highways with minimal latency to destination -- and select these highways by 
imitating other travelers already in the network -- admit an equilibrium point 
provided that the external inflow is bounded above by the min-cut capacity. 
This implies that traffic systems where the users learn through imitation 
transfer the maximum amount of flow that is transferable by a network operating 
at equilibirum. 
This connects our work with classical static flow models used in the 
transportation literature.
Moreover, our results show that when the rate of imitation (namely, the 
frequency at which new users imitate the path preferences of other users) is 
either too small or too large, the equilibrium points may fail to be 
asymptotically stable, thus implying that in unregulated networks the 
congestion state may oscillate around (or escape from) the equilibria.

\textit{Related Work.}
The traffic model proposed here finds its roots in the well-established routing 
game~\cite{TR:05} and corresponding traffic assignment problem~\cite{MP:15}, 
which have been used in the transportation literature to model how travelers 
make decisions in congested traffic.
Recently, this framework has received increased attention with several studies 
investigating the impact of different sources of information on the traffic 
system; e.g., see~\cite{AB-AK-EP-MS:19,AK-AB:20,AF-PG:19,TJ-NL-AB:16} and the 
references therein.
One of the main limitations of this classical approach is that it models 
systems operating at equilibrium, thus neglecting dynamics near these points. 
For this reason, evolutionary dynamics~\cite{JWW:97} have been proposed to study 
the dynamic properties of equilibria~\cite{SF-BV:04,WK-BD-AMB:15}.
Although these works represent a step forward toward understanding the impact 
of advanced navigation systems on traffic patterns, the used models still rely 
on \textit{static descriptions,} where traffic flows propagate instantaneously 
across 
the network. It is immediate to realize that such models are accurate only when 
the routing preferences update at a \textit{slower timescale} than that 
of the traffic dynamics (e.g., when drivers update their path preferences from 
day-to-day as a result of a personal observation)
On the other hand, in modern traffic networks, advanced navigation systems allow 
drivers to update their routing preferences \textit{at the same timescale} as the 
traffic flows, thanks to real-rime traffic state measurements. This connects our 
work with the body of literature on 
dynamic traffic flow models. Our model is a continuous-time version of the 
Cell-Transmission Model~\cite{CFD:95} and related to the model studied 
in~\cite{SC-MA:15}. Dynamic traffic models with static routing preferences 
have been studied in~\cite{GC-KS-DA-MAD-EF:13} using monotonicity, 
in~\cite{SC-MA:16} using mixed monotonicity, in~\cite{GB-FP:20-cdc} using 
passivity. 
Of particular relevance to the framework studied here 
are~\cite{GC-KS-DA-MAD-EF:13,GC-RM:19}. With respect to these works, here we 
study path selection mechanisms governed by the replicator equation and we 
focus on the game-theoretic properties of this model and its stability 
analysis.
This work extends the preliminary work of the authors~\cite{GB-FP:20-cdc} in 
several directions, including a formal proof of uniqueness and evolutionary 
stability of the Nash equilibrium, and a sufficient condition to ensure 
asymptotic stability of the equilibrium point. 
Finally, the recent works~\cite{TT-AK-PF:23,TT-AK-PF:23b} also highlighted 
detrimental effects of navigation systems in a small-scale (two-link) network.

\textit{Contribution.}
The contribution of this work is threefold.
First, we propose a dynamic model derived from the replicator dynamics to 
describe the path selection mechanism underlying drivers' decisions in 
congested traffic. 
We then couple this routing model with a dynamic model of traffic, which 
describes the evolution of traffic flows in the network in relation to the 
instantaneous routing choices. 
Relative to the classical traffic assignment framework, the use of a dynamic 
traffic model describes modern networks where \textit{routing decisions and 
traffic flows update at the same timescale.}
As illustrated in Section~\ref{sec:simulations}, this model allows us to capture 
dynamic behaviors observed in practice, which could not be explained using 
static models~\cite{YN:10,SC-HY-MK:19}.
Second, we study the game-theoretic properties of the equilibria of the 
interconnected model. We show that, under suitable assumptions, an equilibrium 
point exists, is unique, and coincides with an evolutionary stable Nash (or 
Wardrop) equilibrium~\cite{JGW:52}. This relates our work with the 
well-established routing game~\cite{TR:05}.
Third, we study the stability properties of the equilibrium. By using a 
Lyapunov-based reasoning, we derive sufficient conditions under which the 
equilibrium is asymptotically stable. In simulation, we show that the 
conditions are tight and that oscillating trajectories can emerge when our 
conditions do not hold. Intuitively, oscillations originate because 
the population is overreacting to small changes in congestion, more precisely, 
in practice this occurs because individual users update their routing preferences 
without anticipating the preferences of the rest of the population. 
This behavior is consistent with field data (see, e.g., 
\cite{YN:10,SC-HY-MK:19}).

\textit{Organization.}
This paper is organized as follows. Section \ref{sec:2} presents the proposed 
model. Section~\ref{sec:propertiesOfEquilibria} derives conditions for 
existence and uniqueness of an equilibrium point and in 
Section~\ref{sec:asymptotic stability} we study the stability properties of 
the equilibria. 
Section~\ref{sec:simulations} illustrates our findings via numerical 
simulations and Section \ref{sec:conclusion} concludes the paper.

\textit{Notation.} Given  $x\in \real^n, u\in\real^m$, we let 
$(x,u) \in \real^{n+m}$ denote their concatenation; if $n=m, \inprod{x}{u}$ 
denotes the inner product. For symmetric matrix $M,$ $\sbs{\lambda}{max}(M)$ and 
$\sbs{\lambda}{min}(M)$ denote its largest and smallest eigenvalue, 
respectively.

\section{Model of traffic network}
\label{sec:2}

In this section, we present our models of traffic flows and routing 
decisions, and we formalize the problem we study.

\subsection{Traffic flow model}
\label{subsec:traffic_model}
We model a transportation network using a digraph 
$\mc G = (\mc V,\mc E),$ where $\mc V$ is the set of nodes and 
$\mc L$ is the set of links.
In what follows, we let $\mc L = \{1, \dots n\},
n \in \naturalpos.$
For a link $i \in \mc L,$ we denote by $o_i \in \mc V$ its origin node 
and by $d_i \in \mc V$ its destination node.
Motivated by real-world transportation networks with parallel highways, we will 
allow for parallel links, namely, we admit $i,j \in \mc L$ such that $i \neq j$
and have the same origin and destination: $o_i = o_j$  and $d_i=d_j.$
A \textit{path} in $\mc G$ is a sequence of links $\{i_1, i_2, \dots \}$ such 
that the origin node of each link is the destination node of the preceding one. 
Notice that a path may contain repeated links and, going along the path, one 
may reach repeated nodes.
A path is \textit{simple} if it does contain the same link more than once.
The \textit{length} of a path is the number of edges contained in 
$\{i_1, i_2, \dots \}$.
Following the Cell Transmission Model~\cite{CFD:95}, we describe the macroscopic 
behavior of traffic on each link $i \in \mc L$ over time $t \geq 0$ using 
the conservation law: 
\begin{align}
\label{eq:link_dynamics}
\dot x_i(t) &= \sps{f}{in}_i(x(t)) - \sps{f}{out}_i(x(t)),
\end{align}
where $x_i(t)\in \real$ is the traffic density in link $i$, 
$\sps{f}{in}_i(x(t))$ is the traffic inflow entering at upstream, and 
$\sps{f}{out}_i(x_i(t))$ is the traffic outflow exiting at downstream. 
We make the following assumptions on the outflow functions.

\begin{assumption}
\label{as:flowFunctions}
For all $i \in \mc L,$ the outflow function $\sps{f}{out}_i(x)$ depends only on 
the density $x_i,$ namely, $\sps{f}{out}_i(x) = f_i(x_i).$
Moreover, $\map{f_i}{\realnneg}{\realnneg}$ satisfies 
$f_i(x_i)=0$ if and only if $x_i=0,$ it is continuous, and strongly monotone; namely, 
\begin{align}\label{eq:monotonicity_link_flows}
(x_i - \bar x_i) (f_i(x_i) - f_i(\bar x_i)) \geq \mu \vert x_i - \bar x_i\vert^2, 
\end{align}
for some $\mu>0$ and for all $x_i, \bar x_i \in \realnneg$.
\QEDB\end{assumption}
We discuss this assumption in Remark~\ref{rem:increasing_outflows} and we 
illustrate some choices of outflow functions in 
Example~\ref{ex:outflow_functions}. 

Assumption~\ref{as:flowFunctions} guarantees 
that~\eqref{eq:link_dynamics} is a positive system~\cite{LF-SR:00}, namely, for 
every non-negative initial state and every non-negative input at all times, its 
state trajectory is non-negative.
In what follows, for all $i \in \mc L,$ we let 
\begin{align*}
C_i := \sup_{z \in \real} f_i(z),
\end{align*}
and $C = (C_1, \dots, C_n).$
If $f_i(\cdot)$ is unbounded, $C_i = + \infty.$

\begin{remark}[\bf\textit{Validity of Assumption~\ref{as:flowFunctions}}]
\label{rem:increasing_outflows}
It is known (see, e.g., \cite{SC-MA:15}) that the assumption that $f_i(x_i)$ only 
depends on $x_i$ and is strictly increasing is valid provided that we restrict 
our focus to free-flow regimes~\cite{CFD:95}. 
More precisely, it has been shown in~\cite{SC-MA:15} that the free-flow 
equilibrium points of a more complete traffic model (that accounts for congestion 
regimes and backpropagation through the junctions) inherit the same stability 
properties of the model considered here.
Hence, the conclusions drawn here will be applicable also to more complete 
models, provided that their operation is restricted to the free-flow 
regimes~\cite{SC-MA:15}. 
While we acknowledge that accounting for congested regimes is an important 
problem, due to the technical challenges in dealing with a more complete model, 
we leave a generalization of our framework as the focus of future works. 
Regarding the condition $f_i(x_i)=0$ if and only if $x_i=0,$ the ``if'' part 
ensures that no vehicle density can flow out of a link when there is no density 
on it, and the ``only if'' part guarantees that any density is allowed to 
exit.~
\QEDB\end{remark}

\begin{example}[\bf \textit{Flow functions that satisfy Assumption~\ref{as:flowFunctions}}]
\label{ex:outflow_functions}
A class of functions satisfying Assumption~\ref{as:flowFunctions} (and used in, e.g., \cite{GB-FP:19-tits}) is that of linear outflow functions, given by 
\begin{align*}
\sps{f}{out}_i(x_i) = \alpha_i x_i,  & & \alpha_i >0.
\end{align*}
In this case, 
$C_i = + \infty$ and $\mu = \min \{\alpha_i\}_{i \in \mc L}.$
A second class of functions satisfying Assumption~\ref{as:flowFunctions} and 
used in~\cite{GC-KS-DA-MAD-EF:13a}~is 
\begin{align*}
\sps{f}{out}_i(x_i) = C_i (1 - e^{- \beta_i x_i}), && \beta_i>0,
\end{align*}
which is strongly monotone on any bounded set.
\QEDB
\end{example}

Throughout this paper, we will focus on single-commodity networks, 
namely, networks for which there is a single origin node $o$ where 
exogenous traffic flows enter the network, and a single destination node $d,$ 
where flows exit the network. 
We assume that $\mc G$ is outflow-connected, namely, there is a path in 
$\mc G$ from every $i \in \mc L$ to $d.$ To avoid trivial cases, we will also
assume that there exists at least one path from $o$ to $d.$ We denote by 
$\lambda \in \realpos$ the commodity inflow rate at $o$.

To model mass propagation through the nodes, we introduce the scalar routing 
ratios (or routing splits) 
\begin{align*}
\{r_{ij}(t)\}_{i,j \in \mc L}, && t\geq 0,
\end{align*}
where $r_{ij}(t)$ models the fraction of flow exiting link $i$ that proceeds 
toward $j.$
We let $r_{ij}(t)$ be normalized fractions, so that $r_{ij}(t) \in [0,1]$.
Because exchange of flow is allowed only between consecutive links in the 
network, we have $r_{ij}(t)>0$ only if $d_i=o_j.$
Finally, mass is conserved through the nodes when $\sum_ j r_{ij}(t) = 1.$ 
Similarly, we let $r_{oi}(t) \in[0,1]$ be the fraction of exogenous inflow 
$\lambda$ that is routed from the origin node $o$ to link $i;$ analogously, we have $r_{oi}(t)=0$ if $o_i \neq o,$ and $\sum_{i \in \mc L}r_{oi}=1.$
In what follows, it will be useful to combine the network routing ratios into 
a matrix $R(t)=[r_{ij}(t)]\in \real^{n \times n}$ and the routing ratios at the 
origin into a vector $R_o = (r_{o1}, \dots , r_{on}) \in \real^n.$
See Example~\ref{ex:exampleTrafficModel} for an illustration of the model and 
notation.

\begin{remark}[\bf\textit{Temporal dependence in the routing ratios}]
In this discussion, we treated $\{r_{ij}(t)\}_{i,j \in \mc L}$ as time-varying 
quantities; we will see shortly below (cf. Section~\ref{subsec:routing_model}) 
that the time-dependency in $r_{ij}(t)$ implicitly originates as a function of 
the traffic state $x(t).$~
\QEDB\end{remark}

At every node of $\mc G,$ traffic flows are conserved, and thus the inflow to each link
$i \in \mc L$ is given by
$$\sps{f}{in}_i(x) = r_{oi}(t) \lambda +\sum_{j \in \mc L} r_{ji}(t) f_j(x_j(t)).$$ 
By substituting into~\eqref{eq:link_dynamics}, the density on each link evolves 
as:
\begin{align}
\label{eq:linkDynamics}
\dot x_i(t) = r_{oi}(t) \lambda +\sum_{j \in \mc L} r_{ji}(t) f_j(x_j(t)) - f_i(x_i(t)).
\end{align}
By letting 
\begin{align*}
x &:= (x_1, \dots x_n), && 
f(x) := (f_1(x_1),\dots,f_n(x_n))
\end{align*}
be the joint vectors of densities and flows, respectively, the network state 
evolves according to:
\begin{align}
\label{eq:networkDynamics}
\dot x(t) = (R(t)^\tsp-I) f(x(t)) + R_o(t) \lambda.
\end{align}
We illustrate this traffic model in Example~\ref{ex:exampleTrafficModel}. 

\begin{figure}[tb]
\centering
\includegraphics[width=.6\columnwidth]{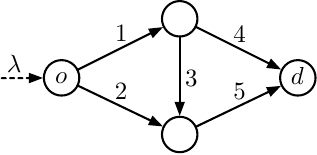}
\caption[]{Graph topology used to illustrate our model. See
examples~\ref{ex:exampleTrafficModel} and \ref{ex:exampleSelectionModel}. 
Nodes labeled by `o' and `d' describe the origin and destination, where exogenous 
inflows enter and exit the network, respectively. The dashed arrow illustrates 
the exogenous inflow.}
\label{fig:wardropNetwork} 
\end{figure}

\begin{example}[\bf \textit{Illustration of traffic flow model}]
\label{ex:exampleTrafficModel}
Consider the network topology in Fig.~\ref{fig:wardropNetwork}. 
The model \eqref{eq:networkDynamics} reads as:
\begin{align*}
\dot x_1 &= -f_1(x_1) + r_{o1} \lambda, & 
\dot x_3 &= -f_3(x_3) + r_{13} f_1(x_1), \\
\dot x_2 &= -f_2(x_2) + r_{o2} \lambda,&
\dot x_4 &= -f_4(x_4) + r_{14} f_1(x_1),\\
\dot x_5 &= -f_5(x_5) + f_2(x_2) + 
f_3(x_3), \hspace{-1cm}
\end{align*}
where time dependencies have been dropped for compactness.
Notice that the routing ratios satisfy $r_{o1}+r_{o2}=1, r_{13}+r_{14}=1.$
\QEDB
\end{example}

%
%\begin{remark}{\bf \textit{(Capturing Backwards Propagation)}}
%\label{rem:backwardsPropagation}
%Notice that our model accounts for saturations in the maximum flow 
%that can be transferred through the highways, but it does not account 
%for spillbacks through the junctions. In this sense, our model is an 
%instance of the Cell Transmission Model~\cite{CFD:95} with infinite 
%supply functions. Notice that this assumption is valid in cases where 
%highways are sufficiently long so that congestion cannot propagate 
%backwards from one highway to another. We leave the investigation of 
%effects due to spillbacks as the topic of future works. 
%%that does not account for spillbacks through the junctions. This 
%%assumption is valid in cases where 
%%
%%
%%More precisely, our model: accounts for limits in the flow that can be 
%%\textit{sent} from a link as described by outflow functions
%%(see Assumption~\ref{as:flowFunctions}), but does account for limits 
%%in the flow that can be \textit{received} by a link since it does not
%%include link supply functions (as done in e.g. \cite{CFD:95}).
%\QEDB
%\end{remark}

\subsection{Congestion-responsive path selection model}
\label{subsec:routing_model}
\begin{figure}[t]
\subfigure[]{\includegraphics[width=.45\columnwidth]{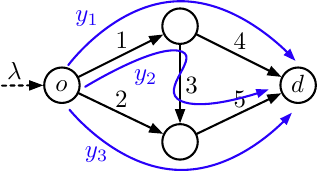}} 
\hfill
\centering \subfigure[]{\includegraphics[width=.45\columnwidth]{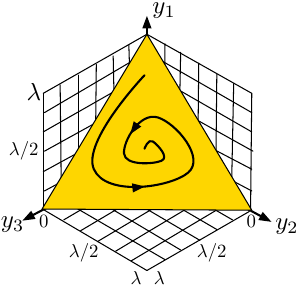}} 
\caption{(a) Demanded path flows: $y_1(t),y_2(t),y_3(t)$ describe the fraction of 
$\lambda$ that is routed through paths, respectively, $p_1, p_2, p_3$ at time 
$t.$ 
(b) The demanded path flows follow a model of path selection 
(cf.~\eqref{eq:replicator}) where the preference for a certain path will increase 
or decrease depending on whether that path has a better or worse travel time 
\textit{compared to the network average}.}
\label{fig:simplex}
\end{figure}

We next propose a model to describe the path selection process followed by 
drivers that seek to minimize their travel time to destination.
Let $\mc P$ denote the set of simple paths from $o$ to $d.$ 
We assume that when a vehicle driver (hereafter called a user) enters the network 
at $o,$ they select a path in $\mc P,$ and they follow this path to destination 
without updating it while traversing the network.
To model this process, we introduce the variables 
$\{y_p(t)\}_{p \in \mc P},$  where $y_p(t)$ denotes the fraction of exogenous 
inflow $\lambda$ that is routed through path $p$ at time $t.$
We stress that $y_p(t)$ models a virtual amount of flow that may never be 
observed in the network: indeed, $y_p(t)$ describes the fraction of $\lambda$ 
that is routed through $p$ at time $t,$ but the actual traffic flows on the links 
of $p$ will be determined by the traffic flow model, as described shortly below.  
Hence, in what follows, we call $y_p(t)$
\textit{demanded path flow} for path $p$ (for a discussion on this wording 
choice, see Remark~\ref{rem:demands_vs_actual} shortly below).
See Fig.~\ref{fig:simplex} for an illustration.
Then, the set of admissible path flow demands is the scaled simplex:
\begin{align}\label{eq:simplex}
\Delta := \{ y \in \realnneg^p: \sum_{p \in \mc P} y_p = \lambda\}.
\end{align}
For link $i \in \mc L,$ we let 
\begin{align}\label{eq:link_flow_demand}
y_i^l(t):= \sum_{p \in \mc P : i \in p} y_p(t),
\end{align}
be the \textit{demanded link flows}. 
Similarly to the demanded path flows, the demanded link flow $y_i^l(t)$ describes 
the fraction of $\lambda$ that is routed through link $i$ at time $t.$
In vector form, $y(t) = (y_1(t), \dots , y_{\vert \mc P \vert}(t))$ 
and $y^l(t) := (y_1^l(t), \dots , y_n^l(t)).$
Notice that~\cite[Thm 2.2]{MP:15} guarantees that for any $y(t) \in \Delta,$
$y^l(t)$ is uniquely determined.

To every link $i \in \mc L,$ we associate a latency function $\ell_i^l (x_i)$ 
mapping traffic density into latency, and describing the travel time or latency  
required to traverse that link.  
With this notation, the total \textit{demanded traffic latency} for path $p$ 
is given by the sum of latencies of the links in that path:
\begin{align}\label{eq:path_latency}
\ell_p(x) := \sum_{i \in p} \ell_i^l(x_i).
\end{align}
In vector form, $\ell(x) := (\ell_1(x), \dots, \ell_{\vert \mc P\vert}(x))$ and 
$\ell^l(x) := (\ell_1^l(x_1), \dots, \ell_{n}^l(x_n)).$
Motivated by~\cite{DB:76}, we make the following assumption on the latency 
functions. 
\begin{assumption}\label{as:continuity_latencies}
For all $i \in \mc L,$ 
$\map{\ell_i^l}{\realnneg}{\realnneg}$ is non-negative, continuous, and such that 
\begin{align}\label{eq:time_tends_infinity}
\lim_{x_i \rightarrow f_i^\inv(C_i)} \ell_i^l(x_i) = + \infty.
\end{align}
%, and  non-decreasing.
\QEDB\end{assumption}

Assumption~\ref{as:continuity_latencies} is very mild, as it requires that 
every link has a non-negative travel time that varies smoothly as a function of 
the traffic densities and that tends to infinity as the link flow approaches the 
flow capacity; we refer to~\cite{DB:76} for a detailed discussion on the validity 
of this assumption.

We consider a model where the vector of path preferences $y(t)$ is continuously 
updated over time based on the traffic state of the network.
We adopt a model of path selection where the preference for a certain path 
$p\in \mc P$ will increase or decrease depending on whether that path has a 
better or worse travel time \textit{compared to the network average}. 
To this end, we model the time-evolution of the flow demands using 
the replicator dynamics~\cite{JH-KS:03}:
\begin{align}
\label{eq:replicator}
\dot y_p(t) = y_p(\bar \ell(x(t),y(t)) - \ell_p(x(t))),
\end{align}
where
\begin{align}\label{eq:replicator_ell_bar}
\bar \ell(x(t),y(t)) = \lambda^\inv \sum_{p \in \mc P} y_p(t)
 \ell_p(x(t)),
\end{align}
is the average latency of traversing the network from $o$ to $d.$
Equation~\eqref{eq:replicator} states that the growth rate of flow demand for 
path $p$ is proportional to the difference between the average latency of 
traversing the 
network $\bar \ell_p(x(t))$ and the latency of that path $\ell_p(x(t)).$ 
We motivate our choice of adopting the replicator dynamics in 
Remark~\ref{rem:motivation_replicator}; we also note that this model 
has been widely adopted in the transportation literature to study dynamics in the 
routing game~\cite{SF-BV:04,BD-WK-AB:14}.

\begin{remark}[\bf\textit{Choice of the replicator dynamics}]
\label{rem:motivation_replicator}
The replicator equation is a deterministic model of imitation, where future 
path preferences are selected by imitating successful path preferences of 
previous users. 
Replicator dynamics have originated in biology, arising in the study of animal 
behavior and evolution, and researchers later proved that this model is also an 
accurate description of processes governed by machine learning algorithms in 
large populations~\cite{TB-RS:97}. 
Interestingly, this is a good model to describe the outcome of machine learning 
processes as its dynamics hinge on historical data and the paradigm of 
imitation (i.e., users observe others' travel times and change their own 
strategies based on these observations). 
%It is worth noting that alternative models of path selections could 
%also be considered, such as the well-established best-response 
%dynamics~\cite{KS:11}, which rely on the paradigm of innovation as 
%opposed to imitation.
Although our analysis is tailored to the replicator model, other selection 
models could be also considered -- such as the best-response 
dynamics~\cite{EH:99}. It is worth noting that the asymptotic 
properties of the trajectories are common across several different models of 
selection: for instance, \cite{EH:99} showed that noisy versions of the best-response dynamics have the same qualitative properties as the replicator 
dynamics. 
\QEDB\end{remark}

\begin{example}[\bf \textit{Illustration of the path selection model}]
\label{ex:exampleSelectionModel}
Consider the network illustrated in Fig.~\ref{fig:motivatingExample} and 
discussed in Example~\ref{ex:exampleTrafficModel}. This graph includes three 
simple paths $\mc P = \{p_1, p_2, p_3 \}$ (see Fig.~\ref{fig:simplex}(a))
given by
\begin{align*}
p_1 = (1,4), &&
p_2 = (1,3,5), &&
p_3 = (2,5).
\end{align*}
The demanded path flows $y_1, y_2,y_3$ are scalar quantities that model the 
fraction of exogenous inflow $\lambda$ that is routed through, respectively, 
paths $p_1,p_2,p_3.$
According to~\eqref{eq:link_flow_demand}, the demanded flows on the links 
$y_1^l, y_2^l,y_3^l, y_4^l, y_5^l$ can be computed from the demanded flows on 
the path as follows:
\begin{align*}
y_1^l = y_1 + y_2, &&
y_2^l = y_3 &&
y_3^l = y_2, &&
y_4^l = y_1, &&
y_5^l = y_2 + y_3.
\end{align*}
In words, the above relationships state that the flow on each link is the sum 
of the flows on paths passing through that link.
The demanded traffic latencies of the paths~\eqref{eq:path_latency} are
\begin{align*}
\ell_1(x) & = \ell_1^l(x_1) + \ell_4^l(x_4), \quad\quad\quad\quad
\ell_3(x) = \ell_2^l(x_2) + \ell_5^l(x_5),\\
\ell_2(x) & = \ell_1^l(x_1) + \ell_3^l(x_3) + \ell_5^l(x_5).
\end{align*}
Namely, the latency of each path is the sum of latencies of all links that 
compose that path. 
The average latency of traversing the network~\eqref{eq:replicator_ell_bar} is:
\begin{align*}
\bar \ell(x,y) = \lambda^\inv(y_1 \ell_1(x) + y_e \ell_2(x) + y_3 \ell_3(x)),
\end{align*}
and models the latency required to traverse the network, averaged over all 
paths. The replicator model~\eqref{eq:replicator} proposed to describe the path 
selection process in this case reads as:
\begin{align*}
\dot y_1 &= y_1 (\bar \ell(x,y) - \ell_1(x)), \quad\quad
\dot y_2 = y_2 (\bar \ell(x,y) - \ell_2(x)),\\
\dot y_3 &= y_3 (\bar \ell(x,y) - \ell_3(x)).
\end{align*}
In words, the preference of users for a certain path grows proportionally to 
the difference between the average delay in the network the delay of that path.
See Fig.~\ref{fig:simplex}(b).
\QEDB\end{example}

It is important to recall some important properties of the replicator 
model~\eqref{eq:replicator} that will be used throughout this paper. 
First, \eqref{eq:replicator} satisfies $\sum_{p \in \mc P} \dot y_p(t)=0$ 
at all times, and thus the simplex $\Delta$ is forward invariant. Namely, if 
$y(0)\in \Delta,$ then, $y(t) \in \Delta$ for all $t>0$.
Second, the boundary faces 
\begin{align*}
\text{bf}_p \Delta := \{y: \sum_{p \in \mc P} y_p = \lambda, y_i=0\}, 
&& p \in \mc P,
\end{align*}
are also forward invariant, and so are the boundary $\text{bd} \Delta$ (i.e., 
the union of all the boundary faces) and the interior $\text{int} \Delta$ (the 
subset satisfying $y_i > 0 ~\forall i$).
It is worth noting that, if $y(0)\in \text{int}\Delta,$ the trajectories 
of~\eqref{eq:replicator} may converge to the boundary only for 
$t \rightarrow +\infty$ and are confined to $\text{int} \Delta$ for all finite 
$t.$

Importantly, these properties imply that if the initial condition $y(0)$ is 
such that $y_p(0)\in \text{bf}_p$ for some $p \in \mc P,$ the replicator 
equation will ignore $y_p$ (namely, $y_p(t)=0$ for all $t \geq 0$). This fact 
implies that one can define a new set of $\vert \mc P \vert-1$ dimensional 
dynamics~\eqref{eq:replicator} where the variable $y_p$ is removed, and 
the trajectories of the $\vert \mc P \vert$ dimensional and 
$\vert \mc P \vert-1$ dimensional dynamics coincide at all times with the 
additional condition $y_p(t)=0$ for all $t \geq 0.$
Motivated by this observation, in what follows it will be convenient to 
restrict the state space $\Delta$ of~\eqref{eq:replicator} to the sub-simplex 
$\Delta'$ given by the support of the vector of initial conditions $y(0):$
\begin{align}\label{eq:simplex_prime}
\Delta':= \{ y \in \realnneg^p: \sum_{p \in \mc P} y_p = \lambda,
y_p=0 ~\forall p: y_p(0)=0 \}.
\end{align}

\subsection{Combined model of traffic with congestion-responsive routing}

In this section, we connect the traffic flow model~\eqref{eq:networkDynamics} 
with the path selection model~\eqref{eq:replicator} to derive a model of 
traffic network with congestion-responsive routing. 
The key observation to relate the two models is that the set of demanded link 
flows $y(t)$ implicitly determines the routing ratios $r_{ij}(t),$ as described 
next.
For a link $j \in \mc L,$ let 
$\theta_j := \sum_{i \in \mc L: o_i = o_j} y_i^l$ denote the total demanded 
flow flowing through its origin node $o_j$. Then, given $y \in \Delta',$ we let 
the routing ratios depend on the demanded traffic flows as follows:
\begin{align}
\label{eq:flow_demand_to_routing}
r_{ij}(y)
= \begin{cases}
0 & \text{ if } o_j \neq d_j,\\
y_j^l/\theta_j & \text{ if } o_j = d_i \text{ and }
\theta_j>0,\\
\frac{1}{\vert \{ k \in \mc L : o_k = o_j\}\vert} & \text{otherwise,}
\end{cases}
\end{align}
where $y^l_j$ is implicitly obtained from $y$ 
using~\eqref{eq:link_flow_demand}.
The model~\eqref{eq:flow_demand_to_routing} states that the outflow exiting 
link $i$ splits among the available downstream links proportionally to the 
total flow demand on each downstream link, provided that each downstream link 
carries a nontrivial amount of flow, and is split uniformly among the 
downstream links otherwise. Notice that other allocation rules may be 
considered (e.g., where splits are non-uniform when $\theta_j=0$).

\begin{figure}[tb]
\centering
\includegraphics[width=\columnwidth]{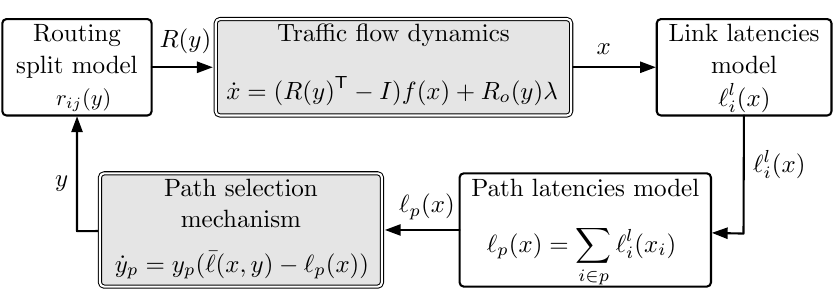} 
\caption[]{The proposed model couples a compartmental-like model of traffic 
flows (see \textit{``Traffic flow physics''}) with an economic model of route 
selection (see \textit{``Path selection mechanism''}). 
Grey-shaded blocks illustrate dynamic models while white blocks illustrate 
algebraic relationships.}
  \label{fig:feedbackInterconnection} 
\end{figure}

By combining \eqref{eq:networkDynamics}, \eqref{eq:path_latency}, 
\eqref{eq:replicator}, and \eqref{eq:flow_demand_to_routing} we obtain the 
following joint traffic flow model with congestion-responsive routing:
\begin{subequations}
\begin{align}
\dot x(t) =  T(x(t), y(t)),\label{eq:interconnectedSystem-a}\\
\dot y(t) = F(x(t),y(t)), \label{eq:interconnectedSystem-b}
\end{align}
\label{eq:interconnectedSystem}
\end{subequations}
\!\!\!\! where $\map{T}{\realnneg^n\times \Delta'}{\real^n}, \map{F}{\realnneg^n\times \Delta'}{\real^p},$ are defined entry-wise, for all 
$i \in \until n$ and $p \in \mc P,$ as:
\begin{align*}
T_i(x,y) &= r_{oi}(y) \lambda +\sum_{j \in \mc L} r_{ji}(y) f_j(x_j) - f_i(x_i),\\
F_p(x,y) &= \eta ~\!y_p(\bar \ell(x,y) - \ell_p(x)). 
\end{align*}
Here, scalar $\eta>0$ is a design parameter that we have introduced to modify 
the rate at which path preferences are updated. 
When the equation~\eqref{eq:interconnectedSystem-b} describes the behavior of 
users following routing recommendations provided by a navigation system, $\eta$
can be modified by deciding the frequency at which travel recommendations are 
updated. For this reason, in what follows, we refer to $\eta$ to as 
\textit{imitation rate.}
We illustrate the interconnection \eqref{eq:interconnectedSystem} and the 
quantities that establish the coupling between the two models in 
Fig.~\ref{fig:feedbackInterconnection}.

We next introduce some basic notation that will be used in the remainder. 
Since~\eqref{eq:flow_demand_to_routing} ensures conservation of flows at the 
nodes, it guarantees that the vector of link flows $y^l$ is an equilibrium 
for~\eqref{eq:networkDynamics}, namely,
\begin{align}\label{eq:equilibrium_demand_flows}
0 = (R(y)^\tsp-I) y^l + R_o(y) \lambda.
\end{align}
From~\eqref{eq:equilibrium_demand_flows}, we deduce that a set of demanded path 
flows $y$ implicitly defines a set of \textit{demanded densities} corresponding 
to these flows. These are defined as:
\begin{align*}
\varphi (y) := f^\inv(y^l), \text{ where }
y^l_i= \sum_{p \in \mc P : i \in p} y_p, ~\forall i \in \mc L,
\end{align*}
and $\map{f^\inv}{\realnneg^n}{\realnneg^n}$ denotes the entrywise inverse 
function of $f(\cdot).$ 
In words, the function $\varphi(y)$ maps a vector of demanded flows into the 
corresponding (demanded) densities.
Similarly to demanded flows, demanded densities are virtual densities, 
which may never be observed in the network, and that model the amount of 
traffic density needed to support the instantaneous demanded flows $y(t).$

We conclude this section by stressing that the flows on the links $f(x(t))$ 
imposed by the traffic dynamics differ from the demanded flows $y^l(t),$ which 
are imposed by the path selection model. We discuss in 
Remark~\ref{rem:flow_demands_vs_actual_flows} the important differences between 
these two quantities.

\begin{remark}[\bf\textit{Demanded flows vs actual flows}]
\label{rem:flow_demands_vs_actual_flows}
\label{rem:demands_vs_actual}
It is important to stress a conceptual difference between ``demanded'' traffic 
variables and traffic variables imposed by the traffic dynamics. 
Regarding traffic flows, the vector of \textit{traffic flows} $f(x(t))$ 
describes the flows on the 
links imposed by the traffic dynamics; on the other hand, the vector of 
\textit{demanded traffic flows} $y^l(t)$ describes the fraction of flow demand 
$\lambda$ entering at $o$ and that is routed to the links based on an economic 
process of path selection.
Analogously, the \textit{traffic densities} $x(t)$ are quantities that are 
imposed by the physics, while the \textit{demanded traffic densities} 
$\varphi(y)$ are virtual quantities describing the densities associated with 
the the traffic demand. 
Importantly, the \textit{traffic latencies} $\ell(x)$ describe the actual travel 
latencies imposed by the physics, which in general differ from the 
\textit{demanded traffic latencies} $\ell(\varphi(y)).$
Notice that the two quantities converge to each other as the dynamics 
of the traffic physics become infinitely fast. 
This discrepancy differentiates our framework from the classical routing 
game~\cite{MP:15}, where the dynamics of traffic  are 
infinitely fast. 
\QEDB\end{remark}

\subsection{Connections with game-theoretic framework}
\label{sec:nash_equilibrium}
In this section, we show that our framework can be related to a population 
game~\cite{WS:15}. This will allows us to connect our setting to the routing 
game~\cite{TR:05} and to relate the equilibirum points of the
model~\eqref{eq:interconnectedSystem} to Wardrop equilibria~\cite{JGW:52}.

The replicator equation~\eqref{eq:replicator} naturally defines an associated 
population game~\cite{WS:15}, as described next. 
A (cost-minimization) population game is defined by the triple 
$(\mc S, \mc X, \kappa),$ where $\mc S$ is a set of pure strategies, $\mc X$ is 
a (generalized) simplex, and $\map{\kappa}{\mc X}{\real^{\vert \mc S\vert}}$ is 
a vector-valued cost function describing the cost associated with each 
strategy, see~\cite[Sec.~13.2]{WS:15}. 
The replicator equation~\eqref{eq:replicator} implicitly defines a population 
game defined by 
\begin{align}\label{eq:game_definition}
\mc S=\mc P, && 
\mc X = \Delta', && 
\kappa(y) = \ell \circ \varphi(y),
\end{align}
which in what follows we denote by 
$\mc R_{\Delta'} := (\mc P, \Delta', \ell \circ \varphi)$.
In line with the existing literature~\cite{WS:15}, we will call a vector of the 
simplex $y \in \Delta'$ a (mixed) \textit{strategy}.
To this end, we will say that a strategy $\sbs{y}{br}\in \Delta'$ is a 
\textit{best reply} to 
$y$ if:
\begin{align*}
\sbs{y}{br}^\tsp ~\ell(\varphi(y)) \leq w^\tsp \ell(\varphi(y)), 
& & \forall w\in \Delta'.
\end{align*}
Associated with $\mc R,$ we have the following classical notion.
% of Nash equilibrium.
\begin{definition}[\bf \textit{Nash Equilibrium}]
\label{def:nash_equilibrium}
A vector $y^*\in \Delta'$ is said to be a \textit{Nash equilibrium} of 
$\mc R_{\Delta'}$ if
\begin{align}\label{eq:nash_equilibirum}
\inprod{y^*}{\ell(\varphi(y^*))} \leq \inprod{y}{\ell(\varphi(y^*))}, & & \forall y \in \Delta'.
\end{align}
\QEDB\end{definition}
In other words, a Nash equilibrium is a best reply to itself.
By noting that $y^\tsp \ell(\varphi(y))$ is the average population latency or 
cost
(cf.~\eqref{eq:replicator_ell_bar}), a Nash equilibrium describes a situation 
where the population has no incentive to deviate away from strategy $y$ as any 
other strategy will yield a non-smaller latency. 
Nash equilibria are used to describe routing games governed by selfish 
individuals, where each individual chooses their path to minimize their travel 
cost.

A very useful reformulation of the notion of Nash equilibrium is that of 
Wardrop equilibrium~\cite{AH-PM:85}: $y$ is a Wardrop equilibrium if, for all 
$p \in \mc P,$ 
\begin{align}\label{eq:wardropcondition}
y_p>0  \text{ implies } \ell_p(\varphi(y)) \leq 
\ell_{p'}(\varphi(y)), ~\forall p' \in \mc P.
\end{align}
In line with the findings of~\cite{AH-PM:85}, in what follows we will use the 
wording Nash equilibrium and Wardrop equilibrium interchangeably.

A desirable property for Nash equilibria is that of 
\textit{evolutionary stability.} 
Intuitively, a strategy $y \in \Delta'$ is evolutionary stable if it is a Nash 
equilibrium and small perturbations from this strategy have a strictly larger 
average latency.

\begin{definition}[\bf \textit{Evolutionary stable point}]
A vector $y\in \Delta'$ is said to be an \textit{evolutionary stable point} of 
$\mc R_{\Delta'}$ if it is a Nash equilibrium and, for all 
$w \in \Delta', w \neq y,$
\begin{align}\label{eq:evolutionary_stability}
w^\tsp \ell(\varphi(y))=y^\tsp \ell(\varphi(y)) 
\text{ implies }
w^\tsp \ell(w) > y^\tsp \ell(\varphi(w)).
\end{align}
\QEDB\end{definition}
In words, $y$ is evolutionary stable if any other best response $w$ to $y$ is 
not a Nash equilibrium. 
It is worth stressing that evolutionary stability is a property of the game 
$\mc R$ as it is defined independently of the choice of the vector 
field in~\eqref{eq:replicator}.

We conclude this section with an important observation, which highlights the 
novelty of the model in Section~\ref{sec:2} with respect to the classical 
routing game framework~\cite{TR:05}. 
We remark that, in the routing game, both the traffic and path selection 
mechanisms operate at the Nash equilibrium~\cite{TR:05} at all times. 
This requirement implicitly makes two highly limiting assumptions: 
(i) the highways have trivial (infinitely fast) dynamics so that the traffic 
flows can be modeled as an algebraic map $\varphi(y)$ of the flow demands; 
(ii) there are no transients in the path selection process, so that the path 
preferences can be described as a Nash equilibrium~\eqref{eq:nash_equilibirum} 
at all times. 
Remarkably, when the framework of the routing game framework was 
derived in the 1950s~\cite{TR:05}, travelers could update their routing 
preferences only from day to day and networks would often 
operate near equilibrium as traffic demands would change slowly. In contrast, 
in modern networks, travelers can update their routing preferences at a fast 
timescale, as they have access to instantaneous real-time traffic information, 
and traffic demands are highly dynamic. Hence, we conjecture that the model 
proposed here is a more accurate description of modern traffic systems.

\section{Properties of the equilibrium points}
\label{sec:propertiesOfEquilibria}
In this section, we study the properties of the equilibrium points 
of the interconnection~\eqref{eq:interconnectedSystem}. 
We begin by noting that solutions to~\eqref{eq:interconnectedSystem} are 
well-defined, as formalized next.

\begin{proposition}[\bf \textit{Well-posedness of solutions}]
\label{prop:existenceUniquenessSolutions}
Let Assumptions~\ref{as:flowFunctions} and \ref{as:continuity_latencies} be 
satisfied, $x(0) \in \realnneg^n$, and $y(0) \in \Delta'$. 
Then, there exists a unique solution $(x(t),y(t)), t\geq 0,$ to the initial 
value problem~\eqref{eq:interconnectedSystem}.
Moreover, $(x(t),y(t)) \in \realnneg^n \times \Delta'$ for all $t\geq 0$.
\QEDB\end{proposition}
\begin{proof}
Existence and uniqueness of the solutions follow from the Lipschitz continuity 
of the vector fields in~\eqref{eq:interconnectedSystem}. The claim 
$x(t) \in \realnneg^n$ follows from Assumption~\ref{as:flowFunctions}, and 
$y(t) \in \Delta'$ follows from $\sum_{p \in \mc P} \dot y_p=0.$
\end{proof}

%The above proposition shows that the solutions to the interconnected 
%dynamics~\eqref{eq:interconnectedSystem} are well-posed. 

\subsection{Existence of fixed points}

We begin by investigating under what conditions the interconnected 
model~\eqref{eq:interconnectedSystem} admits equilibrium points. 
Interestingly, we will show that their existence depends solely on the 
magnitude of external inflows entering the network. To this end, the min-cut 
capacity of the traffic flow model is:
\begin{align*}
\sbs{C}{cut} = \min_{\substack{\mc S \subseteq \mc V: \\ o \in \mc S, d \not \in \mc S}} \sum_{\substack{i \in \mc L:\\ o_i \in \mc S, d_i \not \in \mc S}} C_i.
\end{align*}
Notice that $\sbs{C}{cut}$ may or may not be finite, precisely, 
$\sbs{C}{cut} \in [0, + \infty].$

\begin{proposition}[\bf \textit{Existence of equilibria}]
\label{prop:existence_equilibria}
Let Assumptions~\ref{as:flowFunctions} and \ref{as:continuity_latencies} be 
satisfied.
If $\lambda < \sbs{C}{cut},$ then the interconnected 
system~\eqref{eq:interconnectedSystem} admits an equilibrium point that is a Nash 
equilibrium. 
Conversely, if $\lambda > \sbs{C}{cut},$  then, no equilibirum point exists 
for~\eqref{eq:interconnectedSystem}.
\QEDB\end{proposition}

\begin{proof}
\textit{(Case $\lambda < \sbs{C}{cut}$)}
To prove this implication, we show the existence of a point that satisfies the 
Wardrop conditions and that is an equilibrium 
of~\eqref{eq:interconnectedSystem}. 
Following~\cite[Thm.~2.1]{MP:15}, a vector of path flows 
$\bar y \in \real^{\vert \mc P \vert}$ is a Wardrop 
equilibrium if and only if it satisfies the first-order optimality conditions of 
the following optimization problem:
\begin{subequations}
\begin{align}
\min_{y_1, \dots, y_{\vert \mc P \vert} \in \real} \quad \sum_{i \in \mc L}  \int_0^{y^l_i} & \ell_i^l(s) ds,\label{eq:optimization_traffic_assignment-a}\\
\text{s.~to}~~~\quad\quad \sum_{p \in \mc P} y_p &= \lambda,\label{eq:optimization_traffic_assignment-b}\\
y_p &\geq 0, ~\forall p \in \mc P,\label{eq:optimization_traffic_assignment-c}\\
 \sum_{p \in \mc P: i \in \mc P} y_p &= y^l_i,~\forall i \in \mc L,\label{eq:optimization_traffic_assignment-d}
%y^l_i & < C_i.\label{eq:optimization_traffic_assignment-e}
\end{align}
\label{eq:optimization_traffic_assignment}
\end{subequations}
\hspace{-.3cm} In~\eqref{eq:optimization_traffic_assignment}, $y^l_1, \dots , y^l_n$ are dependent 
variables (describing link flows) that are uniquely determined  
by~\eqref{eq:optimization_traffic_assignment-d} 
(see~\cite[Thm 2.2]{MP:15}). 
Since the objective function is continuous and non-decreasing
(cf.~Assumption~\ref{as:continuity_latencies}), according 
to Weierstrass' Theorem, it admits a minimum provided that the feasible set is 
closed, bounded, and nonempty.
The feasible set of~\eqref{eq:optimization_traffic_assignment} is unbounded, but 
from the positiveness of the latency functions, one may add the constraint 
$y_p \leq d_p,$ where $d_p>0$ is sufficiently large. Hence, the feasible set can be 
made closed and bounded without affecting the solution 
of~\cite[Thm 2.1]{MP:15}. 
Since $\lambda<\sbs{C}{cut},$  by the max-flow min-cut 
theorem~\cite[Thm 4.1]{RKA-TLM-JBO:88}, the feasible set is nonempty. 
Hence, by Weierstrass' Theorem, the game $\mc R$ admits a Nash equilibrium.

Let $y^*$ denote a Nash equilibrium of $\mc R$; we next show that $y^*$ is an 
equilibrium flow for~\eqref{eq:interconnectedSystem-a}. 
Let $R(y^*)$ be the routing matrix obtained from $y^*$  
via~\eqref{eq:flow_demand_to_routing}; by 
using~\eqref{eq:equilibrium_demand_flows}, we have
$$(R(y^*)^\tsp - I) \varphi(y^*) + \lambda =0,$$ 
and thus we conclude that the pair $(x^*, y^*), x^* := \varphi(y^*),$  is an 
equilibrium of~\eqref{eq:interconnectedSystem-a}. 
We are left to show that $(x^*, y^*),$  is also an equilibrium 
of~\eqref{eq:interconnectedSystem-b}.
Since $y^*$ is a Nash equilibrium, it satisfies:
\begin{align*}
\ell_p(\varphi(y^*)) = c, & & \forall p\in \mc P: \bar y_p>0, 
\end{align*}
and thus we have 
$\bar \ell(\varphi(y^*) ,y^*) = c.$ This proves that $(x^*, y^*),$  is an 
equilibrium of~\eqref{eq:interconnectedSystem-b}.

\textit{(Case $\lambda > \sbs{C}{cut}$)}
By contradiction, assume that an equilibrium point $(x^*, y^*)$ exists. 
Because the replicator equation guarantees $y(t)\in \Delta' ~ \forall t\geq 0,$ we must have
\begin{align*}
\sum_{p \in \mc P} y_p^* &= \lambda \text{ and }
 y_p^* \geq 0,  &&\forall p \in \mc P.
\end{align*}
Under these two conditions, the max-flow min-cut theorem is applicable, which 
guarantees that, for some $i \in \mc L,$
\begin{align}
\label{eq:aux_contradiction}
\sum_{p \in \mc P: i \in \mc P} y_p^*  >C_i,
\end{align}
but this contradicts the equilibrium condition
$(R(y^*)^\tsp - I) \varphi(y^*) + \lambda =0,$  thus proving the claim.
\end{proof}

In words, fixed points exist when the external flow demand is 
bounded above by the min-cut capacity; moreover, at least one equilibrium point 
is a Wardrop equilibrium.
This has two important implications. First, it shows that our model is 
consistent with the classical literature, in particular, with the widely 
established notion of Wardrop equilibrium. Importantly, while Wardrop 
equilibria were developed for static models operating at equilibrium, our model 
instead is a dynamic generalization of this classical framework~\cite{MP:15}. 
Second, the result relates our work with the fundamental bound concerning the 
maximum amount of flow transferable by a static graph (as given by the max-flow 
min-cut theorem~\cite{{RKA-TLM-JBO:88}}): it shows that traffic systems where 
users learn through imitation can transfer, asymptotically, the same amount of 
flow as static graphs with free routing. 
This implies that imitation-based selection benefits the overall traffic 
system, enabling it to transfer the maximum admissible amount of flow.
%the fundamental bound characterized by the max-flow min-cut theorem is met. 
We remark that this property is in contrast with dynamic traffic flow models 
with static routing, which may not admit equilibrium points even when 
$\lambda < \sbs{C}{cut}$ (see, e.g., \cite{CFD:95}, 
\cite[Prop.~2]{SC-MA:15}).

\subsection{Conditions for uniqueness of the Nash equilibrium}

While Proposition~\ref{prop:existence_equilibria} guarantees existence of a 
Nash equilibrium, it remains unclear whether such an equilibrium 
is unique or evolutionary stable. We address this aspect next.

\begin{proposition}[\bf \textit{Uniqueness and evolutionary stability}]
\label{prop:nash_evolutionary_stability}
Let Assumptions~\ref{as:flowFunctions}--\ref{as:continuity_latencies} be 
satisfied and $\mc R_{\Delta'}$ be the game induced by~\eqref{eq:replicator} and 
defined as in~\eqref{eq:game_definition}.
Further, assume that the latency functions are strictly monotone, 
namely, for all $i \in \mc L,$
\begin{align}\label{eq:monotonicity_link_latencies}
(x_i - \bar x_i) (\ell^l_i(x_i) - \ell^l_i(\bar x_i)) >0,
\end{align}
for all $x_i, \bar x_i \in [0,C_i), x_i \neq \bar x_i$.
Then, the game $\mc R_{\Delta'}$ admits a unique Nash equilibrium. Moreover, such 
equilibrium is evolutionary stable. 
\end{proposition}

The following lemma is instrumental for the proof. 

\begin{lemma}[\bf \textit{Strict monotonicity of the flow latencies}]
\label{lem:strict_monotonicity_path_flow}
Under the assumptions of Proposition~\ref{prop:nash_evolutionary_stability}, the 
demanded path flow latency functions are strictly monotone, namely, 
\begin{align}\label{eq:monotonicity_flow_latencies}
\inprod{y  - \bar y}{\ell(\varphi(y)) - \ell(\varphi((\bar y))} >0, && \forall y, \bar y \in \Delta', y  \neq \bar y.
\end{align}
\QEDB \end{lemma}

\begin{proof}
We have:
\begin{align*}
\inprod{y  - \bar y&}{\ell(\varphi(y)) - \ell(\varphi((\bar y))} \\
&= \sum_{p \in \mc P} (y_p - \bar y_p)(\ell_p(\varphi(y)) - \ell_p(\varphi(\bar y)))\\
&= \sum_{p \in \mc P} (y_p - \bar y_p) \left( \sum_{i \in p} 
\ell_i^l(\varphi_i(y)) - \sum_{i \in p} \ell_i^l(\varphi_i(\bar y)) \right)\\
&= \sum_{i \in \mc L}  \sum_{p \in \mc P: i \in p} (y_p - \bar y_p) \left( 
\ell_i^l(\varphi_i(y)) -  \ell_i^l(\varphi_i(\bar y)) \right)\\
&= \sum_{i \in \mc L}  (y_i^l - \bar y_i^l) \left( 
\ell_i^l(\varphi_i(y)) -  \ell_i^l(\varphi_i(\bar y)) \right),
\end{align*}
where the second identity follows from~\eqref{eq:path_latency} and the fourth 
identity from~\eqref{eq:link_flow_demand}. 
Next, let $i \in \mc L$ be fixed, and distinguish among three cases. 
(Case 1) Assume $y^l_i>\bar y^l_i,$ we have:
\begin{align}\label{eq:aux_inequality_vz}
(y_i^l - \bar y_i^l) \left( 
\ell_i^l(\varphi_i(y)) -  \ell_i^l(\varphi_i(\bar y)) \right) >0,
\end{align}
since $\phi_i$ and $\ell^l_i$ are strictly increasing. 
(Case 2) Assume $y^l_i<\bar y^l_i.$ In this case, \eqref{eq:aux_inequality_vz} 
also holds since $\varphi_i$ and $\ell^l_i$ are strictly increasing. 
(Case 3) Assume $y^l_i=\bar y^l_i.$ We have:
\begin{align*}
(y_i^l - \bar y_i^l) \left( 
\ell_i^l(\varphi_i(y)) -  \ell_i^l(\varphi_i(\bar y)) \right) =0.
\end{align*}
Since $y \neq \bar y,$ there exists at least one link $i \in \mc L$ for 
which~\eqref{eq:aux_inequality_vz} is satisfied, from which we conclude 
that~\eqref{eq:monotonicity_flow_latencies} holds.
\end{proof}

We are now ready to prove the proposition

\begin{pfof}{Proposition~\ref{prop:nash_evolutionary_stability}}
Let $y^*\in \Delta'$ denote a Nash equilibrium of $\mc R_{\Delta'}$ and 
$y \in \Delta'.$ 
%and for all $i \in \until \mc L$
%\begin{align}\label{eq:link_flow_demand_aux}
%z_i = \sum_{p \in \mc P : i \in p} y_p, & &
%v_i = \sum_{p \in \mc P : i \in p} w_p.
%\end{align}
Using~\eqref{eq:nash_equilibirum}, we have
\begin{align*}
\inprod{y^*}{\ell(\varphi(y^*))}
+ \inprod{y}{\ell(\varphi(y))}  \leq 
\inprod{y}{\ell(\varphi(y^*))} + 
\inprod{y}{\ell (\varphi(y))},
\end{align*}
by re-arranging:
\begin{align}\label{eq:inequality_wy}
\inprod{y}{\ell(\varphi(y))} &\geq \inprod{y^*}{\ell(\varphi(y^*))}
+ \inprod{y}{\ell(\varphi(y))  - \ell(\varphi(y^*)}\nonumber\\
&> \inprod{y^*}{\ell(\varphi(y^*))} 
+ \inprod{y^*}{\ell(\varphi(y))  - \ell(\varphi(y^*)}\nonumber\\
&=  \inprod{y^*}{\ell(\varphi(y))},
\end{align}
where the second row follows from~\eqref{eq:monotonicity_flow_latencies}
and the third row follows from re-arranging the terms.
Inequality~\eqref{eq:inequality_wy} proves~\eqref{eq:evolutionary_stability}, 
thus showing that $y^*$ is evolutionary stable. 
Finally, since the above condition holds for all $y \in \Delta',$ no other point 
in $\Delta'$ other than $y^*$ can satisfy~\eqref{eq:nash_equilibirum}, thus 
proving uniqueness.
\end{pfof}

\begin{figure}[tb]
\centering
\includegraphics[width=.9\columnwidth]{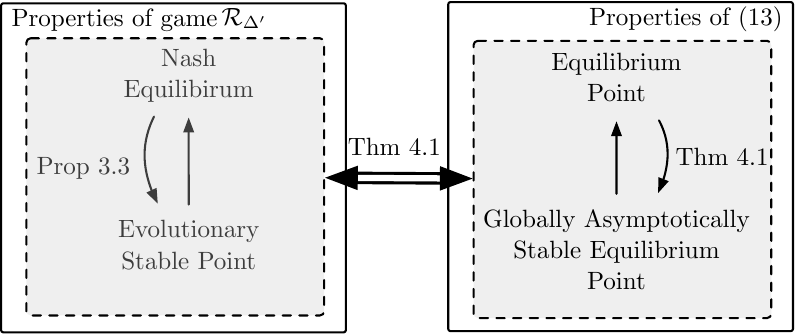}
\caption[]{Summary of the properties established in the main 
results of this paper. Proposition~\ref{prop:nash_evolutionary_stability} shows existence of a unique 
Nash equilibirum and that such equilibirum is evolutionary stable; 
Theorem~\ref{thm:stability_interconnection} shows that the unique 
Nash equilibrium of $\mc R_{\Delta'}$ is also an equilibrium point of~\eqref{eq:interconnectedSystem} and that such a point is globally asymptotically stable.}
\label{fig:relations_results} 
\end{figure}

Proposition~\ref{prop:nash_evolutionary_stability} shows that, under an 
additional monotonicity requirement on the latency functions, the game 
$\mc R_{\Delta'}$
admits a unique Nash equilibrium that is evolutionary stable.
See Fig.~\ref{fig:relations_results}.
We stress that uniqueness and evolutionary stability are properties of the Nash 
equilibrium of the game $\mc R_{\Delta'},$ (and not of the joint 
dynamics~\eqref{eq:interconnectedSystem}).  However, we will show in the next 
section that these properties can be harnessed to study the asymptotic 
properties of the trajectories of~\eqref{eq:interconnectedSystem}.

\section{Asymptotic stability of the Nash equilibrium}
\label{sec:asymptotic stability}
In this section, we study the stability properties of the equilibrium points of 
the interconnection~\eqref{eq:interconnectedSystem}. 
To proceed, we reinforce Assumption~\ref{as:continuity_latencies} as follows.

\begin{assumption}\label{as:strong_monotonicity_latencies}
The conditions in Assumption~\ref{as:continuity_latencies} are satisfied. 
Moreover, the latency functions are strongly monotone, namely, there exists 
$\sigma>0$ such that
\begin{align}\label{eq:strong_monotonicity_link_latency}
(x_i - \bar x_i) (\ell^l_i(x_i) - \ell^l_i(\bar x_i)) \geq 
\sigma \vert x_i - \bar x_i\vert^2,
\end{align}
for all $x_i, x_i' \in [0,C_i)$ and $i \in \mc L.$
\QEDB\end{assumption}

In words, the assumption asks that the latency functions grow at least linearly 
with the traffic densities\footnote{Indeed, strong monotonicity of $\ell^l_i(x_i)$ is equivalent to imposing that 
$\ell^l_i(x_i)-\sigma x_i$ is a monotone function (this follows by rewriting the 
inequality as 
$(x_i - \bar x_i)((\ell^l_i(x_i)-\sigma x_i) - (\ell^l_i(\bar x_i)-\sigma x_i)) \geq 0).$}; the parameter $\sigma$ quantifies the ``steepness'' 
of the density-latency maps. 
In what follows, we will interpret $\sigma$ as a free parameter, which can be 
tuned by a system planner to improve the efficiency of a traffic system 
modeled by~\eqref{eq:interconnectedSystem}.
The following result characterizes the asymptotic behavior 
of~\eqref{eq:interconnectedSystem}.

\begin{theorem}[\bf \textit{Stability of interconnected system}]
\label{thm:stability_interconnection}
Let Assumption~\ref{as:flowFunctions} and \ref{as:strong_monotonicity_latencies}
hold and $\lambda < \sbs{C}{cut}.$ Let $(x(t),y(t))$ denote the solution 
of~\eqref{eq:interconnectedSystem} with initial conditions $(x(0),y(0)),$ 
$\Delta'$ the restricted simplex induced by $y(0),$ $\mc R_{\Delta'}$ the game 
defined by \eqref{eq:game_definition}, and $y^*$ the unique Nash equilibrium of 
$\mc R_{\Delta'}.$
There exists $\sigma^*, \eta_1, \eta_2 >0 $ such that if $\sigma$ and $\eta$ 
satisfy:
\begin{align}\label{eq:condition_sigma_eta}
\sigma > \sigma^*, & & 
\eta \in [\eta_1, \eta_2],
\end{align}
then, for any $(x(0),y(0))\in \realnneg^n \times \Delta',$ 
\begin{align*}
\lim_{t \rightarrow \infty} \norm{(x(t),y(t)) - (x^*,y^*)} = 0,
\end{align*}
where $x^* := \varphi(y^*).$
\QEDB\end{theorem}

The following lemma is a minor extension of 
Lemma~\ref{lem:strict_monotonicity_path_flow} under 
Assumption~\ref{as:strong_monotonicity_latencies}, and is instrumental for the 
proof. 

\begin{lemma}[\bf \textit{Strong monotonicity of the flow latencies}]
\label{lem:strong_monotonicity_flow_latencies}
When Assumptions~\ref{as:flowFunctions} and 
\ref{as:strong_monotonicity_latencies} hold, the path flow latency functions are strongly monotone, namely, there exists $\sigma>0:$
\begin{align}\label{eq:strong_monotonicity_flow_latencies}
\inprod{y  - \bar y}{\ell(\varphi(y)) - \ell(\varphi((\bar y))} \geq \sigma \norm{y - \bar y}^2, && \forall y, \bar y \in \Delta'.
\end{align}
\QEDB \end{lemma}
%The proof of this lemma follows the same steps as 
%Lemma~\ref{lem:strict_monotonicity_path_flow}, with the only difference that 
%inequalities are strengthen using 
%Assumption~\ref{as:strong_monotonicity_latencies}.

\begin{pfof}{Theorem~\ref{thm:stability_interconnection}}
Our proof technique relies on showing that the potential function
\begin{align*}
V(x,y) := V_x(x) + V_y(y),
\end{align*}
where $V_x(x)$ is a potential function for~\eqref{eq:interconnectedSystem-a} and
$V_x(x)$ is a potential function for~\eqref{eq:interconnectedSystem-b} strongly 
decreases along the trajectories of~\eqref{eq:interconnectedSystem} and achieves 
its minimum at $(\varphi(y^*), y^*).$
We will use the following compact notation:
\begin{align*}
A(y):= R(y)^\tsp-I, & & 
\phi(y):= - A(y)^\inv R_o(y)\lambda,
\end{align*}
with $\phi(y)=(\phi_1(y), \dots ,\phi_n(y)).$
Since $\mc G$ is outflow connected, \cite[Thm.~3]{JJ-CS:93} guarantees that 
$A(y)$ is invertible for any $y.$
Since $-A(y^*)$ is a nonsingular 
M-matrix, \cite[Prop.~I${}_{25}$]{RP:77} guarantees the existence of a positive 
diagonal matrix $D = \diag(d_1, \dots, d_n)$ such that
\begin{align}\label{eq:definition_Q}
Q = -(A(y^*) D + DA(y^*)),
\end{align}
is symmetric and positive definite. Let 
\begin{align}\label{eq:Vx}
V_x(x) := 2 \sum_{i \in \mc L} d_i^\inv \int_0^{x_i} f(s) -\phi_i(y^*) ds.
\end{align}
The time-derivative of $V_x(x)$ along the solutions 
of~\eqref{eq:interconnectedSystem} is:
\begin{align}\label{eq:aux_dotVx}
\dot V(x) &= 2 (f(x) - \phi(y^*)) D^\inv (A(y) f(x) + R_o(y)\lambda)\nonumber\\
&= 
2 (f(x) - \phi(y^*))^\tsp D^\inv (A(y^*) f(x) + R_o(y^*)\lambda)\nonumber\\
&\quad\quad + 2 (f(x) - \phi(y^*))^\tsp D^\inv (\psi_x(y) - \psi_x(y^*))\nonumber\\
&= - (f(x) - \phi(y^*))^\tsp Q (f(x) - \phi(y^*))\nonumber\\
&\quad\quad + 2 (f(x) - \phi(y^*))^\tsp D^\inv 
(\psi_x(y) - \psi_x(y^*))\nonumber\\
&\leq - \mu \sbs{\lambda}{min}(Q) \norm{x - \varphi(y^*)}^2 + k \norm{x - \varphi(y^*)}\norm{y-y^*}\nonumber\\
&\leq - \frac{\mu \sbs{\lambda}{min}(Q)}{2} \norm{x - \varphi(y^*)}^2 + \frac{k^2}{2\mu \sbs{\lambda}{min}(Q)} \norm{y-y^*}^2.
\end{align}
Here, in the second row, we used the compact notation
$$\psi_x(y):=A(y) f(x) + R_o(y)\lambda,$$
the third row follows from~\eqref{eq:definition_Q}. 
The fourth row follows by using the Cauchy-Schwarz inequality and by noting that 
$\psi_x(y)$ is Lipschitz continuous in $y,$ uniformly in $x,$ 
%(since $f(x)$ is bounded by Assumption~\ref{as:flowFunctions}), 
and by letting 
$k = 2 \norm{D^\inv}L_\psi L_f,$ where $L_\psi$ and $L_f$ denote the Lipschitz 
constants of $\psi_x(\cdot)$ and $f(\cdot),$ respectively.
The fifth row follows from the inequality $-ax^2+bx \leq b^2/4a$ for 
$a,b>0, x \in \real.$

Next, we let 
\begin{align*}
V_y(y) = \sum_{p \in \mc P} y_p^* \ln\left( \frac{y_p^*}{y_p}\right).
\end{align*}
The time-derivative of $V_y(y)$ along the solutions 
of~\eqref{eq:interconnectedSystem} is given by:
\begin{align}\label{eq:aux_dotVy_1}
\dot V_y(y) &= 
- \eta \sum_{p} y_p^* (\bar \ell(x,y) - \ell_p(x))\nonumber\\
&= - \lambda \bar \ell(x,y) + \eta \sum_{p} y_p^* \ell_p(x)\nonumber\\
&= - \eta \sum_{p} (y_p-y_p^*) \ell_p(x)\nonumber\\
&= - \eta \sum_{p} (y_p-y_p^*) \ell_p(\varphi(y))\nonumber\\
&\quad\quad - \eta \sum_{p} (y_p-y_p^*) (\ell_p(x)-\ell_p(\varphi(y))\nonumber\\
&\leq - \eta \sigma \norm{y-y^*}^2 + \eta L_\ell \norm{y-y^*}\norm{x - \varphi(y)}
\end{align}
Here, the second row follows from $\sum_{p} y_p^*=\lambda;$ the third row follows 
from~\eqref{eq:replicator_ell_bar}; the fourth row from adding and subtracting 
$\ell_p(\varphi(y)).$ The fifth row follows by application of the Cauchy-Schwarz 
inequality, by using continuity of $\ell(\cdot)$ (where $L_\ell$ denotes the 
corresponding Lipschitz constant), and from the following inequality:
\begin{align}\label{eq:aux_strong_bound}
\inprod{y-y^*}{ \ell(\varphi(y))} \geq \sigma \norm{y-y^*}^2.
\end{align}
To prove~\eqref{eq:aux_strong_bound}, since $y^*$ is a Nash equilibrium, we have 
from~\eqref{eq:nash_equilibirum}:
\begin{align*}
\inprod{y^*}{\ell(\varphi(y^*))}
+ \inprod{y}{\ell(\varphi(y))}  \leq 
\inprod{y}{\ell(\varphi(y^*))} + 
\inprod{y}{\ell (\varphi(y))},
\end{align*}
by re-arranging:
\begin{align}\label{eq:inequality_strong}
\inprod{y&}{\ell(\varphi(y))} \geq \inprod{y^*}{\ell(\varphi(y^*))}
+ \inprod{y}{\ell(\varphi(y))  - \ell(\varphi(y^*)}\nonumber\\
&> \inprod{y^*}{\ell(\varphi(y^*))} + c \norm{y-y^*}^2
+ \inprod{y^*}{\ell(\varphi(y))  - \ell(\varphi(y^*)}\nonumber\\
&=  \inprod{y^*}{\ell(\varphi(y))} + c \norm{y-y^*}^2,
\end{align}
where the second row follows from 
Lemma~\ref{lem:strong_monotonicity_flow_latencies}.
This proves~\eqref{eq:aux_strong_bound}.
We can further bound~\eqref{eq:aux_dotVy_1} as:
\begin{align}\label{eq:aux_dotVy_2}
\dot V_y(y) 
&\leq - \eta (\sigma-L_\ell L_\varphi) \norm{y-y^*}^2 \\
&\quad\quad + \eta L_\ell \norm{y-y^*}\norm{x - \varphi(y^*)} \nonumber\\
&\leq - \eta (\frac{\sigma}{2}-L_\ell L_\varphi) \norm{y-y^*}^2 
+ \frac{\eta L_\ell^2}{2 \sigma} \norm{x - \varphi(y^*)}^2, \nonumber
\end{align}
where the first inequality follows from the Cauchy-Schwarz inequality and by 
continuity of $\varphi(\cdot)$ (where $L_\varphi$ denotes the corresponding 
Lipschitz constant), and the second row follows from the inequality 
$-ax^2+bx \leq b^2/4a$ for $a,b>0, x \in \real.$

\begin{figure}[t]
\centering \subfigure[]{\;\;\,\includegraphics[width=.96\columnwidth]{%
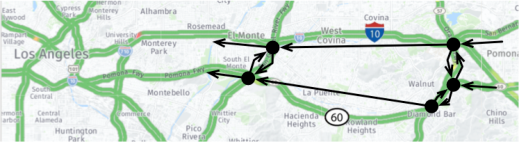}} 
\subfigure[]{\includegraphics[width=\columnwidth]{%
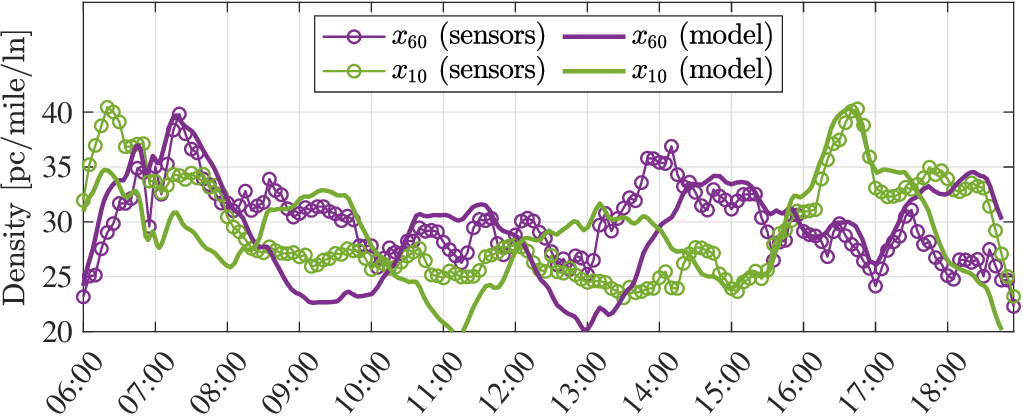}} 
\subfigure[]{\includegraphics[width=\columnwidth]{%
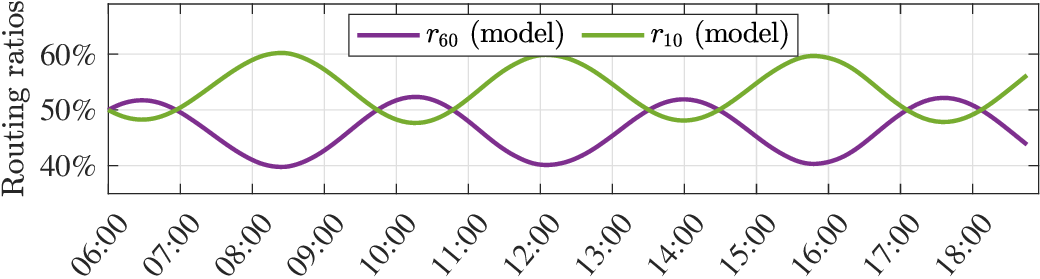}} 
\subfigure[]{\includegraphics[width=\columnwidth]{%
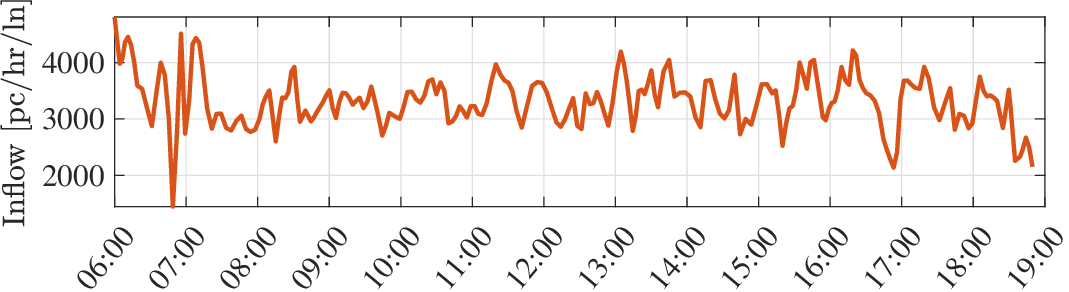}} 
\caption{Time series data for SR60-W and I10-W on March 6, 2020. 
Results obtained by identifying the parameters of~\eqref{eq:interconnectedSystem} 
using a  prediction-correction algorithm that minimizes the fitting error.  
(a) Traffic network and graph. (b) Traffic flow data obtained from sensors 
(continuous lines with circles) and traffic state predicted by our models 
(continuous lines). (c) Routing predicted by our models. (d) Combined traffic 
demand entering at the origin. The data illustrates a case where the trajectories
of~\eqref{eq:interconnectedSystem} oscillate and thus the equilibria lack to be 
asymptotically stable. }
\label{fig:motivatingExample}
\end{figure}

By combining~\eqref{eq:aux_dotVx} and \eqref{eq:aux_dotVy_2} we conclude:
\begin{align}\label{eq:bound_V}
\dot V(x,y) \leq -c_1 \norm{x - \varphi(y^*)}^2 - c_2 \norm{y-y^*}^2,
\end{align}
where the constants $c_1$ and $c_2$ are given by:
\scalebox{.98}{\parbox{.5\linewidth}{
\begin{align*}
c_1 := \frac{\mu \sbs{\lambda}{min}(Q)}{2} - \frac{\eta L_\ell^2}{2 \sigma} , && 
c_2 := \eta\left(\frac{\sigma}{2}-L_\ell L_\varphi\right) - \frac{k^2}{2 \mu \sbs{\lambda}{min}(Q)}.
\end{align*}}}
We thus have that $c_1\geq 0$ and $c_2\geq 0$ when, respectively, 
\scalebox{.97}{\parbox{\linewidth}{
\begin{align}\label{eq:eta1_eta2}
\eta \leq \eta_2:= \frac{\mu \sigma \sbs{\lambda}{min}(Q)}{L_\ell^2}, &&
\eta \geq \eta_1 := \frac{k^2}{\mu \sbs{\lambda}{min}(Q) (\sigma - L_\ell L_\varphi)} .
\end{align}}}
Thus, there exists a feasible choice of $\eta$ that guarantees that $c_1\geq 0$ and $c_2\geq 0$ when $\sigma> \sigma_1 := 2 L_\ell L_\varphi$ and 
\begin{align}\label{eq:sigma_star}
\frac{k^2}{\mu \sbs{\lambda}{min}(Q) (\sigma - L_\ell L_\varphi)} \leq 
\frac{\mu \sigma \sbs{\lambda}{min}(Q)}{L_\ell^2}.
\end{align}
Notice that~\eqref{eq:sigma_star} can always be guaranteed to hold, provided 
that $\sigma$ is chosen sufficiently large.
Altogether this implies that when $\sigma > \sigma^*$ -- where 
$\sigma^* = \max\{\sigma_1, \sigma_2\}$ and $\sigma_2$ is the smallest value of 
$\sigma$ such that~\eqref{eq:sigma_star} holds -- and 
$\eta \in [\eta_1, \eta_2],$ $V(x,y)$ decreases towards its minimum, given by $\xi(x,y)=0,$ which implies 
$(x,y)=(\varphi(y^*),y^*).$ The claim thus follows by application of La Salle's 
invariance principle~\cite[Cor.~4.1]{HKK:96}.
\end{pfof}

\begin{figure}[t]
\subfigure{\includegraphics[width=\columnwidth]{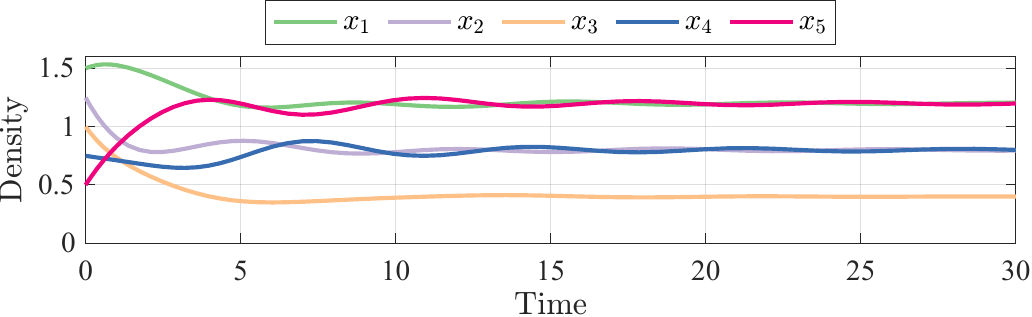}} 
\subfigure{\includegraphics[width=\columnwidth]{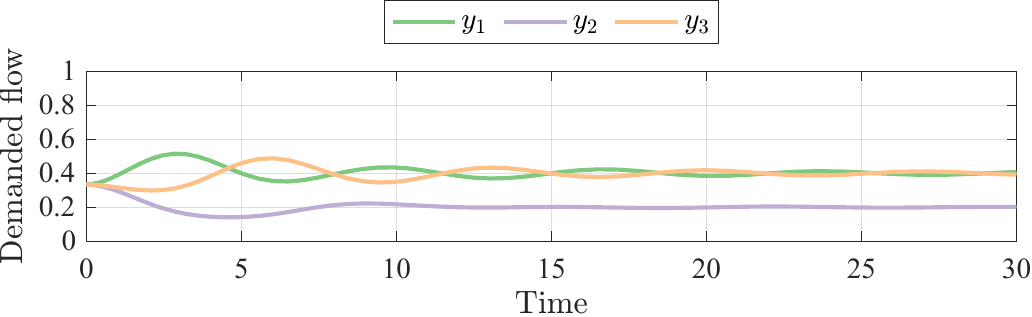}} 
\subfigure{\includegraphics[width=\columnwidth]{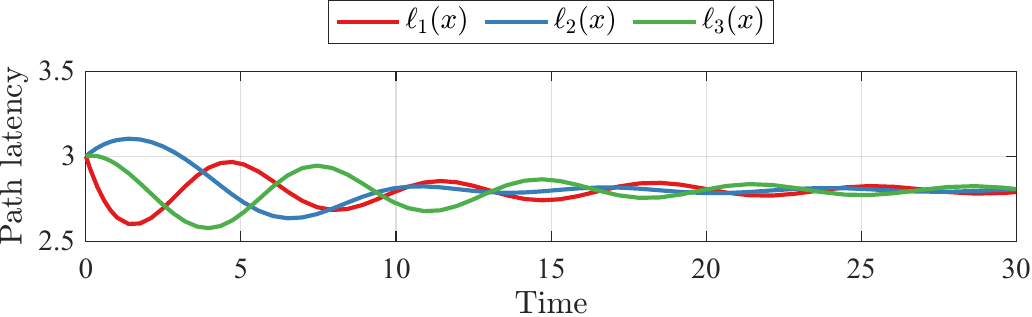}} 
\caption{
Time evolution of the state trajectories of the 
model~\eqref{eq:interconnectedSystem} for the network in 
Fig.~\ref{fig:wardropNetwork} with the choice $\eta=1.$
(Top) Evolution of the traffic density state $x$. (Middle) Evolution of the 
demanded path flow state. (Bottom) Evolution of the travel latencies on the 
paths. The choice $\eta=1$ belongs to the range of stabilizing values 
characterized in Theorem~\ref{thm:stability_interconnection}, and thus 
guarantees that the state asymptotically converges to the Nash equilibrium of 
the underlying game. 
}
\label{fig:braess_eta1}
\end{figure}

We illustrate in Fig.~\ref{fig:relations_results} the relationships between 
implications.
The theorem shows that, provided that the latency functions are sufficiently 
steep and the imitation rate $\eta$ is adequately chosen (as 
in~\eqref{eq:condition_sigma_eta}),  the trajectories 
of~\eqref{eq:interconnectedSystem} converge to the unique Nash equilibrium of 
the game $\mc R_{\Delta'}$ from any initial condition.
We note that, although the statement provides an existence result for 
$\sigma^*, \eta_1, \eta_2,$ an explicit expression for these quantities is 
given in the proof in~\eqref{eq:eta1_eta2} and \eqref{eq:sigma_star}.
Intuitively, \eqref{eq:sigma_star}  states that as $\sigma$ increases, the 
interval $[\eta_1, \eta_2]$ becomes wider since $\eta_1\rightarrow 0$ and 
$\eta_2 \rightarrow+\infty.$ In words, this implies that the steeper the latency 
functions, the more freedom one has in the choice of $\eta.$

Interestingly, the result suggests that asymptotic stability may fail to hold 
when the latency functions are not sufficiently steep, or the imitation rate is 
either too small or too large. Intuitively, when $\sigma$ is small, the 
path selection process is not sufficiently sensitive to variations of traffic 
congestion on the links. On the other hand, when $\eta$ is too large, the 
population is overreacting to small changes in congestion, and individual users 
update their preferences without anticipating the strategy of the rest of the 
population.

\begin{figure}[t]
\subfigure{\includegraphics[width=\columnwidth]{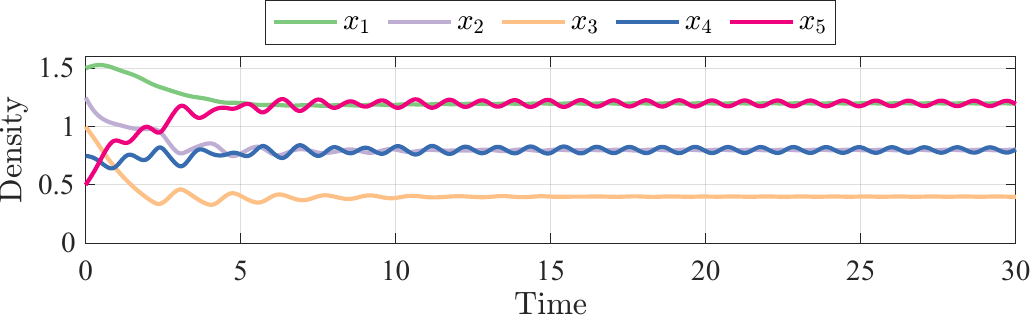}} 
\subfigure{\includegraphics[width=\columnwidth]{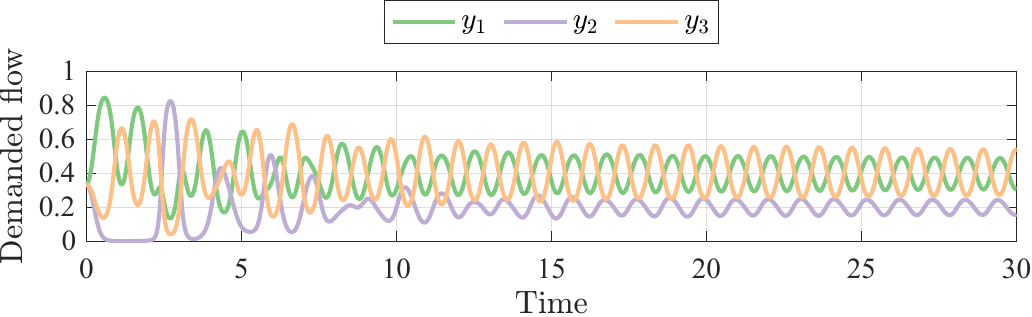}} 
\subfigure{\includegraphics[width=\columnwidth]{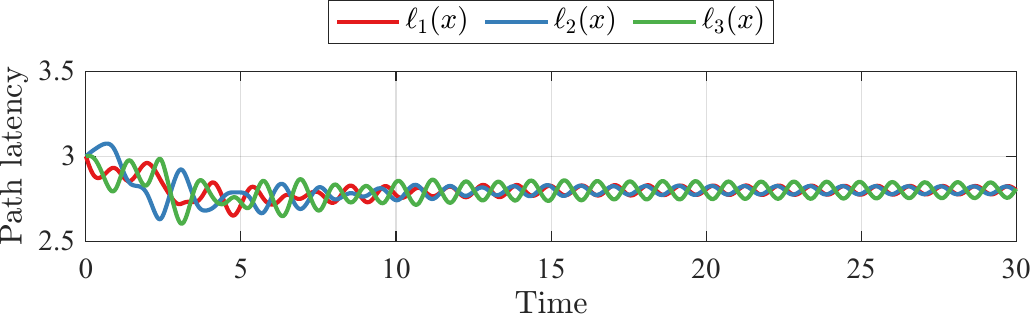}} 
\caption{Time evolution of the state trajectories of the 
model~\eqref{eq:interconnectedSystem} for the network in 
Fig.~\ref{fig:wardropNetwork} with the choice $\eta=1.$
(Top) Evolution of the traffic density state $x$. (Middle) Evolution of the 
demanded path flow state. (Bottom) Evolution of the travel latencies on the 
paths. The choice $\eta=30$ does not belong to the range of stabilizing values 
characterized in Theorem~\ref{thm:stability_interconnection}. As illustrated in 
the simulations, this choice of $\eta$ originates oscillating trajectories, 
describing a condition where users repeatedly switch their path preferences.
}
\label{fig:braess_eta50}
\end{figure}

\section{Simulation Results}
\label{sec:simulations}
This section presents two sets of numerical simulations that illustrate 
our findings.

\subsection{Study case from California SR60-W and I10-W}
Consider the traffic network in Fig.~\ref{fig:motivatingExample}(a), 
which schematizes the west bounds of the freeways SR60-W and 
I10-W in Southern California. 
Let $\sbs{x}{60}$ and $\sbs{x}{10}$ be the average traffic density
in the examined sections of SR60-W  (absolute miles $13.1-22.4$) and 
in the section of I10-W (absolute miles $24.4-36.02$), respectively.
Moreover, let $r_{60}$ (resp. $r_{10}=1-r_{60}$) be the fraction of 
travelers choosing freeway SR60-W over I10-W (resp. choosing 
freeway I10-W over SR60-W) for their commute. 
Fig.~\ref{fig:motivatingExample}(b) illustrates the time-evolution of 
the recorded traffic densities on the two highways on Friday, March 6, 
2020, reconstructed using data from the \textit{Caltrans Freeway 
Performance Measurement System (PeMS)}; in the same figure, we show the 
time-evolution of the state of the interconnected model  
\eqref{eq:interconnectedSystem}. The parameters of the traffic 
system \eqref{eq:networkDynamics} were derived from the nominal 
highway characteristics provided by the PeMS. For the routing model 
\eqref{eq:replicator}, the link latency functions are computed by 
integrating traffic speed data. 
This data illustrates a case where the trajectories 
of~\eqref{eq:interconnectedSystem} oscillate over time, implying that the 
equilibrium points lack to be asymptotically stable; this showcases a scenario 
where the assumptions of Theorem~\ref{thm:stability_interconnection} are not 
satisfied in practice.

\subsection{Illustrative simulations on synthetic model}
Consider the network illustrated in Fig.~\ref{fig:wardropNetwork} and discussed 
in Examples \ref{ex:exampleTrafficModel}-\ref{ex:exampleSelectionModel}. 
Consider a model where $\lambda=1,$ for all $i \in \mc L$ the outflow 
functions are linear $f_i(x_i)=0.5x_i,$ and the latency functions are given by 
$\ell_i(x_i) = x_i, i\in\{1,3,5\}$ and 
$\ell_i(x_i) = 2x_i, i\in\{2,4\}.$
Notice that these choices satisfy Assumption~\ref{as:flowFunctions} 
and~\ref{as:strong_monotonicity_latencies}. 
Proposition~\ref{prop:existence_equilibria} guarantees that the game 
$\mc R_{\Delta'}$ admits an equilibrium point; by 
Proposition~\ref{prop:existenceUniquenessSolutions} such equilibrium is unique 
and evolutionary stable. Solving~\eqref{eq:optimization_traffic_assignment}, 
one obtains the Nash equilibrium $y^*=(2/5, 1/5, 2/5).$
It is then possible to use Theorem~\ref{thm:stability_interconnection} to 
determine values of $\eta$ that guarantee that the trajectories 
of~\eqref{eq:interconnectedSystem} converge to the Nash equilibrium. 
For our choices of functions, one can verify by inspection that 
$\mu = 0.5, L_f=0.5,\sigma=1, L_\ell=2,\sigma=1.$
Moreover, we estimated numerically (sampling each variable uniformly in their 
domain using a Latin Hypercube technique) $L_\varphi=0.125, L_\psi = 1.1547.$ 
We used $D=10^2 I$ and obtained matrix $Q$ (cf.~\eqref{eq:definition_Q}) with 
$\sbs{\lambda}{min}(Q) =20.$ This yields $k=0.2039.$
With these choices, it is easy to see that~\eqref{eq:sigma_star} is verified,
and $\eta_1=2.6667~10^{-6}, \eta_2 = 25.$
Fig.~\ref{fig:braess_eta1} illustrates the state trajectories 
of~\eqref{eq:interconnectedSystem} for $\eta=1.$ As anticipated by 
Theorem~\ref{thm:stability_interconnection}, the state trajectories converge to 
the Nash equilibrium of the game $\mc R_{\Delta'}.$
On the other hand, Fig.~\ref{fig:braess_eta50} illustrates the state 
trajectories of the interconnected system with the choice $\eta=30.$
The simulation demonstrates that an inadequate choice of imitation rate 
$\eta$ leads to trajectories that oscillate over time and not approach the 
Nash equilibrium. The drawbacks of this oscillating phenomenon can be visualized 
by comparing the path latencies illustrated in the bottom figures of 
Fig.~\ref{fig:braess_eta1}  and Fig.~\ref{fig:braess_eta50}.
The choice $\eta=1$ guarantees that all used paths have the same latency at 
equilibrium, thus ensuring that all users experience the same travel time. On 
the other hand, with the choice $\eta=50,$ travel latencies are not homogeneous 
across the three paths, implying that certain users experience a worse travel 
time and higher congestion. From our simulations, we observed that the 
amplitude of oscillating trajectories increases with the flow demand $\lambda,$
thus suggesting that the suboptimality discussed above could deteriorate with 
increased congestion.

\section{Conclusions}\label{sec:conclusion}
This paper proposed a dynamic model of traffic and path selection to describe 
the impact of app-informed travelers in modern traffic networks. We studied the 
properties and stability of the equilibrium points of this model, showing that it 
is consistent with existing studies in transportation. 
Our results suggest that the general adoption of navigation systems enables these 
networks to transfer an amount of flow no smaller than the min-cut capacity, and 
that the equilibrium points are asymptotically stable provided that the latency 
functions are sufficiently sensitive and the imitation rate is adequately chosen.
Future studies should investigate how our conclusions translate to more  general
models that account for bounded supply back propagation through the junctions.
Our results give rise to several opportunities for future work.  By coupling 
these models with common infrastructure control models (such as variable speed 
limits and freeway metering), these results may play an important role in 
designing dynamic controllers for congested infrastructures. Furthermore, our 
models and stability analysis represents a fundamental framework for future 
studies on robustness and security analysis.

\bibliographystyle{myIEEEtran}
%\bibliography{alias,combined,main_GB,GB}
\bibliography{alias,full_GB,GB}

\end{document}